\documentclass[american]{amsart}
\usepackage[T1]{fontenc}
\usepackage[latin1]{inputenc}
\usepackage{amssymb}

\makeatletter

\newcommand{\noun}[1]{\textsc{#1}}

 \theoremstyle{plain}
 \theoremstyle{plain}    
 \newtheorem*{thm*}{Theorem} 
 \theoremstyle{plain}    
 \newtheorem{lem}{Lemma} 
 \theoremstyle{plain}    
 \newtheorem{prop}{Proposition} 
 \theoremstyle{plain}    
 \newtheorem{thm}{Theorem} 

\usepackage[T1]{fontenc}
\usepackage[latin1]{inputenc}

\makeatletter

\usepackage[T1]{fontenc}
\usepackage[latin1]{inputenc}
\usepackage{babel}

\makeatletter

\usepackage[mathscr]{eucal}
\usepackage{amssymb}
\makeatother

\usepackage{babel}
\makeatother

\usepackage{babel}
\makeatother
\begin{document}
   \newcommand{\sA}{{\mathscr{A}}} \newcommand{\sB}{{\mathscr{B}}} \newcommand{\sE}{{\mathscr{E}}} \newcommand{\sF}{{\mathscr{F}}} \newcommand{\sG}{{\mathscr{G}}} \newcommand{\sH}{{\mathscr{H}}} \newcommand{\sJ}{{\mathscr{J}}} \newcommand{\sK}{{\mathscr{K}}} 
\newcommand{\RR}{\mathbb{R}}

\newcommand{\QQ}{\mathbb{Q}}

\newcommand{\EE}{\mathbb{E}}

\newcommand{\Abb}{\mathbb{A}}

\newcommand{\PP}{\mathbb{P}}

\newcommand{\CC}{\mathbb{C}}

\newcommand{\TT}{\mathbb{T}}

\newcommand{\HH}{\mathbb{H}}

\newcommand{\NN}{\mathbb{N}}

\newcommand{\ZZ}{\mathbb{Z}}

\newcommand{\sP}{\mathcal{P}}

\newcommand{\esp}[1]{\EE\left[ #1 \right] }

\newcommand{\pr}[1]{\PP\left[ #1 \right] }
 
\newcommand{\Rp}{\RR_{+}^{*}}
 
\newcommand{\RxR}{\Rp\times\RR}

\newcommand{\affR}{\mathrm{Aff}(\RR)}

\newcommand{\affT}{\mathrm{Aff}(\TT)}

\newcommand{\affO}[1]{\mathrm{Aff}(#1 )}
 
\newcommand{\affq}{\mathrm{Aff}(\QQ_{p})}

\newcommand{\affQ}{\mathrm{Aff}(\QQ)}

\newcommand{\hor}{\mathrm{Hor}(\TT)}

\newcommand{\bT}{\partial^{*}\TT}
 
\newcommand{\bG}{\partial^{*}\Gamma}

\newcommand{\horg}{\mathrm{Hor}(\Gamma)}
 
\newcommand{\horO}[1]{\mathrm{Hor}(#1 )}

\newcommand{\toinf}{\rightarrow\infty}
 
\newcommand{\inv}{^{-1}}
 
\newcommand{\Ind}[1]{1_{#1 }}
 
\newcommand{\nor}[1]{\left| #1 \right| }

\newcommand{\norp}[1]{\left| #1 \right| _{p}}
 
\newcommand{\sumzi}[1]{\sum_{#1 =0}^{\infty}}
 
\newcommand{\dnor}[1]{\left\Vert #1 \right\Vert }

\newcommand{\tm}{\frac{3}{2}}

\newcommand{\stard}{\stackrel{\cdot}{*}}

\newcommand{\supp}{\mathrm{supp}}

\newcommand{\cf}{\wedge}

\newcommand{\e}{\mathrm{e}}

\newcommand{\ver}{\rightarrow}

\newcommand{\Qp}{\QQ_{p}}

\newcommand{\nup}[1]{\nu_{p}\left( #1 \right) }

\newcommand{\nad}[1]{\left\langle #1 \right\rangle }

\newcommand{\sumP}{\sum_{p\in\sP}}

\newcommand{\card}[1]{\textrm{card}\left\{  #1 \right\}  }

\newcommand{\ZP}{\ZZ(P)}

\newcommand{\affP}{\affO{P}}

\newcommand{\gP}{\left( P\right) }

\newcommand{\norP}[1]{\left[ #1 \right] _{P}}

\newcommand{\nX}[1]{\left\Vert #1 \right\Vert }

\newcommand{\bnd}{\mathbf{bnd}}

\newcommand{\bP}{\overline{\sP}}

\newcommand{\pd}{\phi_{p}}

\newcommand{\Pm}{P^{*}}

\newcommand{\naP}[1]{\left\langle #1\right\rangle _{P}^{+}}

\begin{flushright}\today{}\end{flushright}

\title{The Poisson boundary of random rational affinities}

\maketitle
\smallskip{}
\begin{center}\noun{Sara Brofferio}%
\footnote{\noindent Institut für Mathematik C -- Technische Universität Graz\\
 Steyergasse 30 -- A-8010 Graz\\
 brofferio@finanz.math.tu-graz.ac.at\\
http://finanz.math.tu-graz.ac.at/\textasciitilde{}brofferio/

Supported by Marie Curie Fellowship HPMF-CT-2002-02137%
}\end{center}
\smallskip{}

\begin{abstract}
The group of affine transformations with rational coefficients $\affQ$
acts naturally on the real line $\RR$, but also on the $p$-adic
fields $\Qp$. The aim of this note is to show that, for random walks
whose laws have a finite first moment, all these actions are necessary
and sufficient to describe the Poisson boundary, which is in fact
the product of all the fields that contract in mean. 
\end{abstract}

\section*{Introduction}

Random walks are processes on a group $G$ defined as iterated products
of independent and identically distributed random elements and are
a natural probabilistic way to explore the algebraic structures and
their underlying geometry. The complex interaction between these mathematical
objects can be illustrated by the Poisson boundary. The latter can
be defined pure measure theoretically as the space that contains all
the informations on the long range behavior of the random walk, but
it is also the maximal one among the $\mu$-boundaries, which are
the topological $G$-spaces that are stable and contracting under
the action of the random walk. Furthermore it has an interpretation
from an analytic viewpoint, since it provides an integral representation
of all harmonic bounded functions.

The study of random walks on the group of rational affinities $\affO{\QQ}$,
which is the group of transformations of the form $x\mapsto ax+b$
(or equivalently of the matrices $\left[{{a\: b\atop 0\:1}}\right]$)
where the coefficients $a\neq0$ and \textbf{$b$} are rational numbers,
is a good example of how a quite elementary probabilistic process
is related to sophisticated arithmetic spaces. This countable group
has a natural action on the real line $\RR$ and is a dense subgroup
the group of real affine transformations. One can obtain interesting
results concerning the behavior of the random walks on $\affO{\QQ}$
using the powerful theory developed on Lie groups, when no continuity
hypothesis on the measure is assumed (for instance \cite{Ke73}, \cite{BP92},
\cite{BBE} or \cite{Br03b}). Nevertheless the Poisson boundary for
random walk on $\affQ$ can not be studied in such a way. In fact,
while, for random walk with a spread out law on $\affR$, the boundary
is either $\RR$ or trivial (see L.Elie \cite{El84}), V.Kaimanovich
\cite{Ka91} showed that for random walk supported by the group of
affine transformations with dyadic coefficients the Poisson boundary
is either the real line or the other possible completion of the dyadic
line, namely, the 2-adic field $\QQ_{2}$. It was suggested that a
complete understanding of the asymptotic behavior of the random walks
on $\affQ$ could be obtained by considering simultaneously the actions
on $\RR$ and on all $p$-adic fields $\Qp$where $p$ is in $\sP$,
the set of all prime numbers. 

Since the formal structure of the real and $p$-adic fields is similar,
they can often be approached in a similar way and, in order to unify
the notation, it is common to associate the real setting to the {}``prime
number'' $p=\infty$, thus $\QQ_{\infty}=\RR$. Under first moment
conditions, the parameter that determines whether the action of $\affQ$
is contracting in mean on the field $\Qp$ is the $p$-drift \[
\phi_{p}=\int_{\affQ}\ln\norp{a}d\mu(a,b),\]
 where $\mu$ is the step law of the random walk. When $\pd$ is negative,
there is a unique $\mu$-invariant probability measure on $\Qp$,
which is in fact a non-trivial $\mu$-boundary.

The aim of this note is to show that, for all measures with a first
moment on $\affQ$, the Poisson boundary is the product of all $p$-adic
fields with negative drift. We prove the following :

\begin{thm*}
Let $\mu$ be a probability measure on $\affQ$ that is not supported
by an Abelian subgroup and such that\[
\int_{\affQ}\left(\sum_{p\in\sP}\nor{\ln\norp{a}}+\sum_{p\in\sP\cup\left\{ \infty\right\} }\ln^{+}\norp{b}\right)\: d\mu(a,b)<+\infty.\]
 Then there exists a unique $\mu$-invariant probability measure $\nu^{*}$
on the topological product\[
B^{*}=\prod_{p\in\sP\cup\left\{ \infty\right\} :\phi_{p}<0}\Qp,\]
 and the measure space $(B^{*},\nu^{*})$ is the Poisson boundary
of the random walk of law $\mu$. 
\end{thm*}
Furthermore, the measure $\nu^{*}$ carries no point mass except in
the case when $B^{*}$ collapses to a single point, namely when $\pd\geq0$
for all $p\in\sP\cup\left\{ \infty\right\} $. Since the $p$-drifts
have to satisfies to $\phi_{\infty}=-\sum_{p\in\sP}\phi_{p}$, we
deduce that the Poisson boundary is trivial if and only if all $p$-drifts
are null. 

This paper is organized as follows. In section 1, we quickly introduce
the basic concepts of $\mu$-boundary and of Poisson boundary. In
section 2, we summarize well known results on the contracting action
of $\affQ$ on the fields $\Qp$. We deduce that $B^{*}$ is a $\mu$-boundary
that is a good candidate to be maximal. We also observe that, even
though the topological space $B^{*}$ is not locally compact, the
measure $\nu^{*}$ is supported by a set $B_{\mathbf{r}}^{*}$ that
is in fact a restricted topological product of the $\Qp$ with respect
to some of their compact discs, and thus it can be endowed with a
locally compact topology homeomorphic to a sub-space of the Adele
ring. To prove that $(B^{*},\nu^{*})$ is in fact the Poisson boundary,
we use the techniques based on the estimation of the entropy introduced
by Kaimanovich and Vershik \cite{KV83} and Derrienic \cite{Der80},
and, in particular, the criterion on the entropy of the conditional
expectation due to Kaimanovich \cite{Ka00}. Our main tool is the
construction of a suitable family of gauges in terms of what we shall
call an \emph{adelic length} on $\affQ$, based on the arithmetic
\emph{height} of the Adeles (section 3). This permits to estimate
the growth of the random walk and prove some laws of large numbers
(section 4). In section 5, using the projection of $B^{*}$ onto finite-dimensional
$\mu$-boundaries, we can prove that its conditional entropy is zero
and, thus, that it is the Poisson boundary.

In our previous note \cite{Br04}, we studied the Poisson boundary
for measures $\mu$ that are supported by finitely generated subgroups
of $\affQ$, using the Strip approximation criterion \cite{Ka00}.
This technique cannot be applied directly in the present more general
context, since the random walk grows faster and it is not straightforward
to exhibit a global geometrical approximation. On the other hand,
the adelic length provides a suitable tool to control the entropy
of the $\mu$-boundaries and, since the technical arguments turn up
to be quite light, it is likely that this approach can be adapted
to more general algebraic groups over rational numbers.

\section{$\mu$-boundaries and Poisson boundary }

Let $\mu$ be a probability measure on a countable group $G$ and
let $\left\{ g_{n}\right\} _{n}$be a sequence of independent random
elements with law $\mu$ on $G$. Consider the (right) \emph{random
walk} $\left\{ x_{n}\right\} _{n}$ starting at the identity, which
is the process on $G$ defined by \[
x_{n}=g_{1}\cdots g_{n}\qquad\forall n\in\NN.\]
We denote by $(G^{\NN},\PP)$ the probability space of trajectories
of the random walk and by $\EE$ the associated expectation.

Let $B$ a locally compact $G$-space endowed with a $\mu$-invariant
probability measure $\nu$ such that $\PP$-almost surely $x_{n}\nu$
converges vaguely to a Dirac measure (where for every $g\in G$ the
measure $g\nu$ is given by $g\nu(f)=\int_{B}f(gz)d\nu(z)$ ). According
to Furstenberg \cite{Fur73}, the space $(B,\nu)$ is a \emph{$\mu$-boundary}
and the \emph{Poisson boundary} is the maximal of such spaces, namely
it is a $\mu$-boundary such that any other $\mu$-boundary is one
of its measurable $G$-equinvariant quotients.

One can define a measurable map $\bnd=\bnd_{B}$ from the space $(G^{\NN},\PP)$
to the $\mu$-boundary $(B,\nu)$ that associates to a path $\mathbf{x}=\left\{ x_{n}\right\} $
the point $\bnd(\mathbf{x})$ of $B$ such that\[
\lim_{n\toinf}x_{n}\nu=\delta_{\bnd(\mathbf{x})}\qquad\textrm{almost surely}.\]
In other words, the action of the random walk on $B$ contracts to
$\bnd(\mathbf{x})$, which contains all the the informations on the
asymptotic behavior of $x_{n}$ acting on $B.$

As a measure space, the Poisson boundary is unique and there exist
several equivalent constructions for a generic countable group. For
instance, it can be identified with the quotient of the probability
space $(G^{\NN},\PP)$ by the equivalence relation\[
\left\{ x_{n}\right\} _{n}\sim\left\{ x'_{n}\right\} _{n}\Longleftrightarrow\exists k,h\in\NN:x_{n+k}=x'_{n+h}\forall n\in\NN,\]
namely with the measure space that contains all possible long term
behaviors of the random walk. A classical question is to give a tangible
description of this measure space and to recognize when a given topological
(or measure) space, which is known to be $\mu$-boundary, is in fact
the Poisson boundary.

\section{$\mu$-boundaries of $\affQ$}

The group of rational affinities \[
\affO{\QQ}=\left\{ (a,b):x\mapsto ax+b\mid a\in\QQ^{*},\, b\in\QQ\right\} \]
 has by definition an action on the group of rational numbers. However,
$\QQ$ endowed with the discrete topology cannot be a $\mu$-boundary,
because it cannot support a stationary probability measure, except
in degenerate cases.

The action of $\affQ$ on the rational numbers extends naturally to
the real line $\RR$, but also to the $p$-adic numbers $\Qp$ for
all prime numbers $p$. These fields are the completion of $\QQ$
with respect to \emph{$p$-adic norm} \[
\left|q\right|_{p}=p^{-v_{p}(q)},\]
where the \emph{$p$-adic valuation} of an integer $r$ is $v_{p}(r)=\max\left\{ k\in\NN\mid p^{-k}r\in\NN\right\} $
and $v_{p}(r/s)=v_{p}(r)-v_{p}(s),$ while $\norp{0}=0$. The real
and the $p$-adic norms are known to be the only possible norms on
$\QQ$ adapted to its field structure. 

Since the real and $p$-adic fields are formally similar, it is useful
to associate the {}``prime number'' $p=\infty$ to the real setting;
thus $\QQ_{\infty}$ is $\RR$, the Euclidean norm is $\nor{\,\cdot\,}_{\infty}$
and so on. We denote by $\sP$ the set of all true prime numbers and
write $\overline{\sP}=\sP\cup\left\{ \infty\right\} $.

Let us consider a probability measure $\mu$ on $\affQ$ and the associated
random walk $x_{n}$ obtained as the product of the sequence $\left\{ g_{n}=(a_{n},b_{n})\right\} _{n}$
of random affinities with law $\mu$. A simple calculation shows that
\[
x_{n}=(A_{n},Z_{n})=(a_{1}\cdots a_{n},\sum_{k=1}^{n}a_{1}\cdots a_{k-1}b_{k}).\]
 We always suppose that the law $\mu$ is \emph{non-degenerate}, that
is :\[
\pr{a_{1}=1}\not=1\quad\textrm{and }\quad\pr{a_{1}z+b_{1}=z}\not=1\quad\forall z\in\QQ.\]
In fact, whenever this does not hold, the random walk degenerates
either to a sum of independent random variables in $\QQ$ or to a
product of independent elements in $\QQ^{*}$(using the map $(a,b)\mapsto(a,az+b)$).
In both cases the support of $\mu$ generates an Abelian group, and
it is well known that the Poisson boundary is trivial. 

If for some $p\in\overline{\sP}$ the measure $\mu$ has a (logarithmic)
\emph{first $p$-moment}, that is\[
\esp{\ln\norp{a_{1}}+\ln^{+}\norp{b_{1}}}<\infty,\]
 the parameter that determines whether the action on the respective
field $\Qp$ is contracting is the \emph{$p$-drift} \[
\pd=\esp{\ln\norp{a_{1}}}.\]
 In fact one has the following classical results

\begin{lem}
Suppose that $\mu$ is non-degenerate and has a first $p$-moment.

\textbf{a.} If $\phi_{p}<0$, the infinite sum \begin{equation}
Z_{\infty}^{p}=\sum_{k=1}^{\infty}a_{1}\cdots a_{k-1}b_{k}\label{eq-Z-in-Qp}\end{equation}
converges almost surely in $\Qp$ to a random element with law $\nu_{p}$,
which carries no point mass. Furthermore $(\Qp,\nu_{p})$ is a $\mu$-boundary. 

\textbf{b.} If $\phi_{p}\geq0$, there exists no stationary probability
measure on $\Qp$.

\end{lem}
\begin{proof}
For the convenience of the reader, we give a sketch of the proof.

First observe that, since the measure $\mu$ is supposed to be non-degenerate,
no stationary probability measure $\nu$ can carry a point mass. In
fact, suppose that some point of $\Qp$ carries a non-null mass. Let
$M$ be the maximum of such masses and $S=\left\{ z\in\Qp|\nu(\left\{ z\right\} )=M\right\} $.
Then $g\cdot S=S$ for all $g$ in the support of $\mu$. Let $s\in S$,
since the measure is not degenerate there exists $g\in\supp\mu$ such
that $g\cdot s\not=s$. But since each affinity fixes just one point
of $\Qp$, the orbit $\left\{ g^{n}\cdot s\right\} _{n\in\NN}$ is
infinite and thus $S$ should be infinite too, which is absurd. 

\textbf{a.} $\pd<0$. For more details on the real case see \cite{Ver79}
and on the ultra-metric case see \cite{CKW}.

Observe that by the Law of large numbers, the process \[
\norp{a_{1}\cdots a_{n}}=\exp\left(\sum_{i=1}^{n}\ln\norp{a_{i}}\right)\]
converges almost surely to zero with exponential speed (roughly as
$\exp(n\phi_{p})$). On the other hand, since $\ln^{+}\norp{b_{1}}$
is integrable, $\left.\ln^{+}\norp{b_{n}}\right/n$ converges almost
surely to zero. Thus the infinite sum (\ref{eq-Z-in-Qp}) converges,
because its general term goes to zero exponentially. 

Furthermore, $\PP$-almost surely for all $z\in\Qp$\[
x_{n}\cdot z=A_{n}x+Z_{n}\,\rightarrow\, Z_{\infty}^{p}\qquad\textrm{in }\Qp.\]
Thus, by dominated convergence, for every continuous bounded function
$f$ on $\Qp$ \[
x_{n}\nu_{p}(f)=\int_{\Qp}f\left(A_{n}x+Z_{n}\right)\nu_{p}(dx)\rightarrow\delta_{Z_{\infty}^{p}}(f)\qquad\PP-\textrm{almost surely }.\]

\textbf{b.} $\pd\geq0$. Bougerol and Picard \cite{BP92} obtained
an analogous result for stationary sequences of multidimensional real
affinities. We translate here their proof to the case of a sequence
of independent $p$-adic affinities.

Let $p\in\sP$ be a true prime number and suppose that there exists
a stationary probability measure $\nu$ on $\Qp$. Let $f$ be a non-negative
bounded function on $\Qp$ with compact support and consider the process
\[
W_{n}=\int_{\Qp}f(x_{n}\cdot z)d\nu(z).\]
This is a bounded martingale and, thus, it converges almost surely
and in $L^{1}$ to a non-negative random variable $W_{\infty}$. Furthermore,
if $f\neq0$, the random variable $W_{\infty}$ is not null, because
its mean is $\nu(f)$. Let now \[
m_{n}=\max\left\{ \norp{A_{n}},\norp{Z_{n}}\right\} .\]
 Observe that $m_{n}$ is a power of $p$ and that $\norp{m_{n}}=m_{n}^{-1}$.Thus,
the sequences $\left\{ m_{n}A_{n}\right\} _{n}$ and $\left\{ m_{n}Z_{n}\right\} _{n}$
are bounded in $\Qp$. Since $\phi_{p}\geq0$ (whence $\norp{A_{n}}$
is unbounded), there exists a sub-sequence $\left\{ n_{i}\right\} _{i}$
such that $m_{n_{i}}$ diverges to $+\infty$ and such that \[
m_{n_{i}}A_{n_{i}}\ver A\quad\textrm{and }\quad m_{n_{i}}Z_{n_{i}}\ver Z\quad\textrm{in }\Qp,\]
for some $A,Z\in\Qp$. Then, for all $z\not=-Z/A$\[
\lim_{i\toinf}\norp{x_{n_{i}}\cdot z}=\lim_{i\toinf}m_{n_{i}}\norp{m_{n_{i}}A_{n_{i}}z+m_{n_{i}}Z_{n_{i}}}=+\infty\]
 Thus, $\PP$-almost surely\[
W_{\infty}=\lim_{i\toinf}\int_{\Qp}f(x_{n_{i}}\cdot z)d\nu(z)=\lim_{i\toinf}\int_{\Qp}f(x_{n_{i}}\cdot z)\Ind{\left[z=-Z/A\right]}d\nu(z)=0,\]
since $\nu$ has no point mass. Thus we obtained a contradiction.
\end{proof}
It follows from this last lemma that an exhaustive $\mu$-boundary
of $\affQ$ should involve all $p$-adic fields with negative drift.
Let \[
\Pm=\left\{ p\in\bP\mid\mu\textrm{ has a first $p$-moment and }\phi_{p}<0\right\} \]
 and consider the topological sum $B^{*}=\prod_{p\in\Pm}\Qp$ with
the topology $\mathcal{T}_{s}$ generated by the open sets\[
\prod_{p\in S}O_{p}\prod_{p\not\in S}\QQ_{p}\]
where $S\subseteq\Pm$ is finite and the $O_{p}\subseteq\Qp$ are
open. It is easily checked that the action of the random walk on $(B^{*},\mathcal{T}_{s})$
is contracting. In fact, $\PP$-almost surely \begin{equation}
x_{n}\cdot\mathbf{z}=(x_{n}\cdot z_{p})_{p}\ver\mathbf{Z_{\infty}^{*}}=(Z_{\infty}^{p})_{p}\quad\textrm{in }(B^{*},\mathcal{T}_{s})\label{eq:conv-prod}\end{equation}
for all $\mathbf{z}=(z_{p})_{p}\in B^{*}$. Let $\nu^{*}$ denote
the law of $\mathbf{Z_{\infty}^{*}}$.

Observe that whenever $\Pm$ is infinite the space $(B^{*},\mathcal{T}_{p})$
is not locally compact. However it is possible to construct a smaller
locally compact topological $\affQ$-space that supports the measure
$\nu^{*}$.

Since the random variables $Z_{\infty}^{p}$ are almost surely finite,
there exists a sequence $\mathbf{r}=(r_{p})_{p\in\Pm}$ of real positive
numbers greater or equal to $1$ such that\[
\sum_{p\in\Pm}\pr{\norp{Z_{\infty}^{p}}>r_{p}}<\infty.\]
 Thus by the Borel-Cantelli Lemma \[
\pr{\norp{Z_{\infty}^{p}}>r_{p}\textrm{ for an infinite number of $p\in\Pm$}}=0,\]
and the random variable $\mathbf{Z_{\infty}^{*}}$ is almost surely
in the set \[
B_{\mathbf{r}}^{*}=\left\{ \mathbf{z}\in\prod_{p\in\Pm}\Qp:\norp{z_{p}}\leq r_{p}\textrm{ for all $p$ but a finite number}\right\} .\]
This set is locally compact (second countable), if considered as the
restricted topological product of the $\left(\Qp\right)_{p\in\Pm}$
with respect to the discs $D_{p}(\mathbf{r})$ of center $0$ and
radius $r_{p}$ in $\Qp$, that is endowed with the topology $\mathcal{T}_{\mathbf{r}}$
generated by the open sets \[
\prod_{p\in S}O_{p}\prod_{p\not\in S}D_{p}(\mathbf{r})\]
where $S\subseteq\Pm$ is finite and the $O_{p}$ are open subsets
of $\Qp$. 

The topology $\mathcal{T}_{\mathbf{r}}$ is finer than the restriction
of the product topology $\mathcal{T}_{s}$ to $B_{\mathbf{r}}^{*}$,
but the sigma-algebras they generate coincide. Thus $(B_{\mathbf{r}}^{*},\nu^{*})$
and $(B^{*},\nu^{*})$ are the same probability space and, even if
the action on $B_{\mathbf{r}}^{*}$ is not strongly contracting as
in (\ref{eq:conv-prod}), we have the following

\begin{prop}
$(B^{*},\nu^{*})=(B_{\mathbf{r}}^{*},\nu^{*})$ is a $\mu$-boundary.
\end{prop}
\begin{proof}
Let $\overline{B_{\mathbf{r}}^{*}}$ be the one point compactification
of $B_{\mathbf{r}}^{*}$. Thus the sequence $\left\{ x_{n}\nu^{*}\right\} _{n}$
of probability measures on $\overline{B_{\mathbf{r}}^{*}}$ is relatively
compact. Let $\nu'$ be an accumulation point. By (\ref{eq:conv-prod}),
for every bounded function $f$ continuous with respect to $\mathcal{T}_{s}$,
the sequence $x_{n}\nu^{*}(f)$ converges to $f(\mathbf{Z}_{\infty}^{*})$.
Thus $\nu'=\delta_{\mathbf{Z}_{\infty}^{*}}$ on the sigma-algebra
generated by $\mathcal{T}_{s}$ that coincides with the sigma-algebra
generated by $\mathcal{T}_{\mathbf{r}}$.
\end{proof}
\textbf{Remarks.} 1. Since the real and $p$-adic norms satisfy the
following relation \[
\nor{q}_{\infty}=\prod_{p\in\sP}\norp{q}^{-1}\qquad\forall q\in\QQ^{*},\]
if all the $p$-drifts exist, one has \[
\phi_{\infty}=-\sum_{p\in\sP}\phi_{p}.\]
Thus the $\phi_{p}$ cannot be all simultaneously negative and $P^{*}\not=\bP$.
This implies that the space $B^{*}$ and $B_{\mathbf{r}}^{*}$ do
not involve all the possible completions of the rationals numbers.
The Strong approximation theorem (see for instance Cassels \cite{CF}
, page 67) ensures then that diagonal embedding of $\QQ$ in $B_{\mathbf{r}}^{*}$
is always dense. It follows that when the support of the measure $\mu$
generates $\affQ$ as a semi-group, then the support of measure $\nu^{*}$
is the whole of $B_{\mathbf{r}}^{*}$, and thus this boundary is in
some sense minimal.

2. The \emph{Adele ring} $\Abb$ is the restricted topological product
of all $\Qp$ with $p\in\bP$ with respect the disc $D_{p}(1)$ of
center $0$ and radius $1$ (see for instance \cite{CF}, page 63).
Since it is possible to choose the $r_{p}$ in the form $p^{k}$ (thus
$\norp{r_{p}}=r_{p}^{-1}$), the map\[
\begin{array}{ccc}
B_{\mathbf{r}}^{*} & \rightarrow & B_{\mathbf{1}}^{*}\\
(z_{p})_{p} & \mapsto & (r_{p}z_{p})\end{array}\]
is a homeomorphism, which embeds $B_{\mathbf{r}}^{*}$ into a sub-space
of the Adele ring. However, this map is not an isomorphism of $\affQ$-spaces
and, although it is possible to formalize explicit conditions under
which $\mathbf{Z}_{\infty}^{*}$ is almost surely in $\Abb$, this
is not true in a general setting.

\section{Gauges on $\affQ$}

In the previous section we have provided what seems a good candidate
to be the Poisson boundary. To prove that this boundary contains all
the informations on the tail of the trajectories, we need to estimate
the growth of the random walk with respect to the geometry adapted
to this boundary. We have seen that $\affQ$ can be embedded in each
of the affine group over the $p$-adic or real fields, whence it is
natural to use the associated norms. However, since $\QQ$ is dense
in $\Qp$, in order to obtain a topological space that respects the
discrete structure of $\affQ$, one has to consider the diagonal embedding
\begin{eqnarray*}
\affQ & \hookrightarrow & H:=\QQ^{*}\times\Abb\\
(a,b) & \mapsto & (a,(b)_{p}).\end{eqnarray*}
 In fact, since $\QQ$ is discrete in the Adele ring (see \cite{CF})
and $\QQ^{*}$ is endowed with the discrete topology, $\affQ$ is
discrete in $H$. We would like to observe that, since the boundary
$B^{*}$ is not contained in the Adeles, one may be tempted to use
instead of $\Abb$ the restricted topological product of the $\Qp$
with respect to some bigger discs $D_{p}(\mathbf{r})$, but the resulting
embedding would not be discrete.

The space $H$ can be endowed with a group structure by extending
the product on $\affQ$, that is, setting\[
(a,(z_{p})_{p})(a',(z'_{p})_{p})=(aa',(az'_{p}+z_{p})_{p}).\]

For every $q\in\QQ^{*}$ set \[
\nad{q}:=\sum_{p\in\sP}\nor{\ln\norp{q}}.\]
Observe that even if this function may appear exotic, it can be easily
calculated since for every irreducible fraction $\frac{r}{s}$ of
integers, one has $\nad{\frac{r}{s}}=\ln r+\ln s$. 

For all $\mathbf{z}=\left(z_{p}\right)_{p}\in\Abb$, also set

\[
\nad{\mathbf{z}}^{+}:=\sum_{p\in\bP}\ln^{+}\norp{z_{p}}\]
This function,well known in number theory (see for instance Lang \cite{Lan66}),
is called \emph{height}. 

Finally for all $(a,b)\in H$ we define the \emph{adelic length}\[
\nX{(a,b)}=\nad{a}+\nad{b}^{+},\]
which plays, in some way, the role of the word length in this non-finitely-generated
context. 

The function $\dnor{\cdot}$ is not sub-additive, but we have the
following relation\[
\dnor{y_{1}y_{2}}\leq\ln2+2\dnor{y_{1}}+\dnor{y_{2}}\qquad\forall y_{1},y_{2}\in H\]
In fact \[
\nad{a_{1}a_{2}}\leq\nad{a_{1}}+\nad{a_{2}}\]
 and\begin{eqnarray*}
\nad{b_{1}+a_{1}b_{2}}^{+} & = & \ln^{+}\nor{b_{1}+a_{1}b_{2}}_{\infty}+\sum_{p\in\sP}\ln^{+}\norp{b_{1}+a_{1}b_{2}}\\
 & \leq & \ln2+\ln^{+}\left|b_{1}\right|_{\infty}+\ln^{+}\nor{a_{1}b_{2}}_{\infty}+\sum_{p\in\sP}\max\{\ln^{+}\norp{b_{1}},\ln^{+}\nor{a_{1}b_{2}}_{p}\}\\
 & \leq & \ln2+\nad{b_{1}}^{+}+\nad{a_{1}}+\nad{b_{2}}^{+}.\end{eqnarray*}
 Define the gauge $\mathcal{G}^{y}=\left\{ \mathcal{G}_{k}^{y}\right\} _{k\in\NN}$
of center $y\in H$ by setting\[
\mathcal{G}_{k}^{y}=\left\{ g\in\affQ|\nX{g^{-1}y}\leq k\right\} .\]
The sets $\mathcal{G}_{k}^{y}$ are not empty and they exhaust the
whole group. Furthermore their growth is controlled by the following
:

\begin{lem}
The family of gauges $\left\{ \mathcal{G}^{y}\right\} _{y\in H}$
has uniform exponential growth, that is, there exists $C>0$ such
that $\card{\mathcal{G}_{k}^{y}}\leq\e^{Ck}$ for all $y\in H$ and
all $k\in\NN$.
\end{lem}
\begin{proof}
First observe that $\card{q\in\QQ^{*}|\nad{q}\leq k}\leq2\e^{2k}$,
since we remarked that $\nad{\frac{r}{s}}=\ln r+\ln s$ when $r,s\in\NN$.
Also observe that if $q=\frac{r}{s}\in\QQ^{*}$ then\[
\nad{\frac{r}{s}}^{+}=\nad{-\frac{r}{s}}^{+}=\ln s+\left(\ln r-\ln s\right)^{+}\]
thus \[
\frac{\nad{q}}{2}\leq\nad{q}^{+}\leq\nad{q}.\]
 We can easily conclude that \[
\card{\mathcal{G}_{k}^{(1,0)}}=\card{(a,b)\in\QQ^{*}\times\QQ:\nad{a}+\nad{b}^{+}\leq k}\leq2\e^{2k}(2\e^{2k}+1).\]

Take now a generic $y=(a,\mathbf{z})\in H$ . It is known (see for
instance Cassels \cite{CF}, page 65), that since $a^{-1}\mathbf{z}$
is in the Adele ring, there exists $b\in\QQ$ such that $\norp{a^{-1}z_{p}-b}\leq1$
for all $p\in\overline{\sP}$. Let $y'=(a,ab)$ and $t=(1,\left(b-a^{-1}z_{p}\right)_{p})$.
Then \[
\dnor{g^{-1}y'}=\dnor{g^{-1}yt}\leq\ln2+2\dnor{g^{-1}y}+\dnor{t}=\ln2+2\dnor{g^{-1}y},\]
 and thus $\mathcal{G}_{k}^{y}\subseteq\mathcal{G}_{\ln2+2k}^{y'}$.
Finally, since $y'\in\affQ$ and \[
\mathcal{G}_{k}^{y}\subseteq\mathcal{G}_{\ln2+2k}^{y'}=y'\mathcal{G}_{\ln2+2k}^{(1,0)},\]
 the lemma follows.
\end{proof}

\section{Laws of large numbers}

Suppose now that the law $\mu$ of the random walk has a first moment
with respect to the gauge function $\dnor{\cdot}$, that is\[
\esp{\dnor{(a_{1},b_{1})}}=\sumP\esp{\nor{\ln\norp{a_{1}}}}+\sum_{p\in\bP}\esp{\ln^{+}\norp{b_{1}}}<\infty.\]
Observe that this global moment condition implies that $\mu$ has
all first $p$-moments, thus that all $p$-drifts exist, and \[
\sum_{p\in\bP}\left|\pd\right|<\infty.\]
However, this condition is not very strong, since it is equivalent
to ask that the numerators and denominators of $a_{1}$ and $b_{1}$
have finite logarithmic moment.

We are going to control the growth of the random walk \[
x_{n}=(a_{1},b_{1})\cdots(a_{n},b_{n})=(A_{n},B_{n})\]
with respect to $\dnor{\cdot}$ by providing a sequence of points
in $H$ depending only on the boundary point $\bnd(\mathbf{x})$ that
well approximates the path $\mathbf{x}=\left\{ x_{n}\right\} $.

We have already observed that the $p$-adic norm of the linear part
$A_{n}$ of the random walk is just the exponential of a sum of i.i.d.
random variables whose mean is the $p$-drift. Thus we can approximate
$A_{n}$ with the rational number \[
q_{n}:=\prod_{p\in\sP}p^{-[n\frac{\phi_{p}}{\ln p}]},\]
whose $p$-norm is of the order of $\e^{n\pd}$ (where $[x]$ is the
integer part of $x$ and $q_{n}$ is in $\QQ$, as $\pd/\ln p$ converges
to zero when $p$ grow). This approximation holds not only locally
on each field $\QQ_{p}$, but also globally when we consider all fields
together.

\begin{lem}
\label{Fct-LLN-QQ*}$\left.\nad{A_{n}^{-1}q_{n}}\right/n$ converges
in $L^{1}$ to zero.
\end{lem}
\begin{proof}
Observe, for every $p\in\overline{\sP}$ by the ergodic theorem \[
\frac{\ln\norp{A_{n}^{-1}q_{n}}}{n}=\frac{-\sum_{k=1}^{n}\ln\norp{a_{k}}+\phi_{p}-\phi_{p}+\ln p[n\frac{\phi_{p}}{\ln p}]}{n}\rightarrow0\]
 in $L^{1}$. Thus, by dominated convergence, the sequence \[
\esp{\frac{\nad{A_{n}^{-1}q_{n}}}{n}}=\sum_{p\in\sP}\frac{\esp{\nor{\ln\norp{A_{n}^{-1}q_{n}}}}}{n}\]
 converges to zero, because each term of the infinite sum converges
to zero and is dominated by $\esp{\nor{\ln\norp{a_{1}}}}+\left|\phi_{p}\right|$,
which is summable over $p\in\sP$.
\end{proof}
In the first part of this paper, we have shown that the action of
the random walk on the different fields $\Qp$ depends on the sign
of the $p$-drift. Thus, we are going to decompose the Adeles in different
parts and, for all $P\subseteq\bP$, we define a partial height \[
\nad{\mathbf{z}}_{P}^{+}:=\sum_{p\in P}\ln^{+}\norp{z_{p}}\qquad\mathbf{z}\in\prod_{p\in P}\Qp.\]
Suppose that for all $p\in P$ the $p$-drift is negative. Then the
translation component $Z_{n}$ of the random walk converges in $\prod_{p\in P}\Qp$,
which is in fact a $\mu$-boundary. Let then \[
\mathbf{Z}_{\infty}^{P}:=\left(Z_{\infty}^{p}\right)_{p\in P}=\lim_{n\rightarrow\infty}Z_{n}\qquad\textrm{in }\prod_{p\in P}\Qp.\]
Observe that, if $\mathbf{0}:=\left(0\right)_{p}\in\prod_{p\in P}\QQ_{p}$,
then \[
x_{n}^{-1}\cdot\mathbf{Z}_{\infty}^{P}=(g_{1}\cdots g_{n})^{-1}\lim_{k\to\infty}g_{1}\cdots g_{k}\cdot\mathbf{0}=\lim_{k\to\infty}g_{n+1}\cdots g_{k}\cdot\mathbf{0}\cong\mathbf{Z}_{\infty}^{P}\]
where the last equality is in law. Thus we have proved the following
:

\begin{lem}
\label{fct-LGN-Zn-centr}Suppose that $\phi_{p}<0$ for all $p\in P\subseteq\overline{\sP}$
. The sequence\[
x_{n}^{-1}\cdot\mathbf{Z}_{\infty}^{P}\]
 is stationary. Thus, if $\naP{\mathbf{Z}_{\infty}^{P}}$ is almost
surely finite, then $\left.\naP{x_{n}^{-1}\cdot\mathbf{Z}_{\infty}^{P}}\right/n$
converges in probability to zero.
\end{lem}
We would like to remark that, in the general case, this result does
not apply directly to the boundary point in the most complete boundary
$B^{*}$. In fact when $\mathbf{Z}_{\infty}^{*}$ is not contained
in the subspace $B_{\mathbf{1}}^{*}$ of the Adeles, its partial height
is almost surely infinite. We will deal with this problem by projecting
$B^{*}$ on products of finitely many $\Qp$.

To estimate $Z_{n}$ on the other directions we use the fowling 

\begin{lem}
For all $P\subseteq\overline{\sP}$ \[
\pr{\frac{\naP{Z_{n}}}{n}\leq\sum_{p\in P}\phi_{p}^{+}+\varepsilon}\rightarrow1\]
 for all $\varepsilon>0$. 
\end{lem}
\begin{proof}
First observe that if $p\in\sP$, by the ultra-metric property \[
\ln^{+}\norp{Z_{n}}\leq\max_{1\leq k\leq n}\ln^{+}\norp{a_{1}\cdots a_{k-1}b_{k}}=:M_{n}^{p}.\]
 For $p=\infty$, keeping the same notation, one has \[
\ln^{+}\left|Z_{n}\right|_{\infty}\leq\ln n+M_{n}^{\infty}.\]
Let $S_{n}^{p}=\sum_{k=1}^{n}\ln\norp{a_{k}}$. The sequence \[
u_{n}^{p}=\frac{M_{n}^{p}}{n}-\phi_{p}^{+}=\max_{1\leq k\leq n}\frac{\left(S_{k}^{p}+\ln\norp{b_{k}}\right)^{+}}{n}-\phi_{p}^{+}\]
 converges to zero almost surely to zero, because $n^{-1}\left(S_{n}^{p}+\ln\norp{b_{n}}\right)^{+}$
converges to $\phi_{p}^{+}$. Since the $u_{n}^{p}$ is bounded by
$\left|\phi_{p}\right|+n^{-1}\sum_{k=1}^{n}\left(\left|\ln\norp{a_{k}}\right|+\ln^{+}\norp{b_{k}}\right)$,
which converges in $L^{1}$, it is uniformly integrable and, thus
converges to zero also in $L^{1}$. Therefore, since \[
\esp{u_{n}^{p}}\leq\esp{\left|\ln\norp{a_{1}}\right|+\ln^{+}\norp{b_{1}}}+\left|\phi_{p}\right|\textrm{ }\]
and \[
\sum_{p\in P}\esp{\left|\ln\norp{a_{1}}\right|+\ln^{+}\norp{b_{1}}}+\left|\phi_{p}\right|\leq3\esp{\dnor{(a_{1},b_{1})}}<\infty,\]
 the sequence $\left\{ \sum_{p\in P}u_{n}^{p}\right\} _{n}$ converges
to zero in $L^{1}$, when $n\to\infty$. Finally, as \[
\nad{Z_{n}}_{P}^{+}\leq\ln n+\sum_{p\in P}M_{n}^{p}\leq\ln n+n\cdot\left(\sum_{p\in P}u_{n}^{p}+\sum_{p\in P}\phi_{p}^{+}\right),\]
 the lemma follows.
\end{proof}
It is now possible to estimate the growth of the random walk on $\affQ$.

\begin{prop}
\label{pro:LLN-affQ} Let $P$ be a finite subset of $\overline{\sP}$
such that $\phi_{p}<0$ for $p\in P$ and let\begin{eqnarray*}
\pi_{n}=\pi_{n}^{P}:\prod_{p\in P}\Qp & \longrightarrow & H\\
\mathbf{z} & \mapsto & (q_{n},\mathbf{z}\cup\left(0\right)_{p\not\in P}).\end{eqnarray*}
 Then\[
\pr{\frac{\nX{x_{n}^{-1}\pi_{n}(\mathbf{Z}_{\infty}^{P})}}{n}\leq\sum_{p\in P^{c}}\phi_{p}^{-}+\varepsilon}\rightarrow1\]
 for all $\varepsilon>0$.
\end{prop}
\begin{proof}
Observe that\[
\dnor{x_{n}^{-1}\pi_{n}(\mathbf{Z}_{\infty}^{P})}=\nad{A_{n}^{-1}q_{n}}+\nad{x_{n}^{-1}\cdot\mathbf{Z}_{\infty}^{P}}_{P}^{+}+\nad{x_{n}^{-1}\cdot0_{P^{c}}}_{P^{c}}^{+}.\]
For any fixed time $n$, the product $x_{n}^{-1}=g_{n}^{-1}\cdots g_{1}^{-1}$
has the same distribution of the random walk $\check{x}_{n}=(\check{A}_{n},\check{Z}_{n})$
associated with the measure $\check{\mu}$, image of $\mu$ by the
inversion in the group. Thus, since $x_{n}^{-1}\cdot0_{P^{c}}$ has
the same law as $\check{Z}_{n}$ and the $p$-drift associated to
$\check{\mu}$ is $\check{\phi}_{p}=-\phi_{p}$, we can apply the
previous lemmas in order to conclude. 
\end{proof}

\section{The Poisson boundary of $\affQ$}

As announced, we are going to prove that the $\mu$-boundary $B^{*}$
is in fact the Poisson boundary by using the criterion based on the
entropy of the conditional expectation developed by Kaimanovich \cite{Ka00}.
Suppose that the measure $\mu$ has finite entropy \[
-\sum_{q\in G}\mu(g)\ln\mu(g)<\infty.\]
 Consider the family $\PP^{z}$ of probability measures obtained conditioning
measure $\PP$ with respect to the events $\bnd_{B}(\mathbf{x})=z$
and let $\PP_{n}^{z}$ be the corresponding measure on the group,
obtained by the projection $\mathbf{x}\mapsto x_{n}$. Then Theorem
4.6 in \cite{Ka00} says that the $\mu$-boundary $(B,\nu)$ is in
the Poisson boundary if and only if for $\nu$-almost all $z\in B$\[
-\frac{1}{n}\ln\PP_{n}^{z}(x_{n})\ver0\qquad\PP^{z}(d\mathbf{x})-\textrm{almost surely}.\]

\begin{thm}
Suppose that $\mu$ has a first moment with respect to $\dnor{\cdot}$.
Then $(B^{*},\nu^{*})$ is the Poisson boundary.
\end{thm}
\begin{proof}
Since $\mu$ has a first moment with respect to a gauge with exponential
growth, it has finite entropy by \cite{De86}.

Observe that if $\mathbf{Z}_{\infty}^{*}$ were in the Adele ring,
then Proposition \ref{pro:LLN-affQ} would hold also for $P=P^{*}$
and $\left.\dnor{x_{n}^{-1}\pi_{n}(\mathbf{Z}_{\infty}^{*})}\right/n$
would converge to 0 in probability, since $\sum_{p\not\in P^{*}}\pd^{-}=0$.
Then by Theorem 5.4 in \cite{Ka00}, we could directly prove that
$(B^{*},\nu^{*})$ is the Poisson boundary. But since $\nad{\mathbf{Z}_{\infty}^{*}}_{P^{*}}^{+}$
is not necessarily finite, we need to be more careful.

Let $P$ be a finite subset of $P^{*}$. For $\mathbf{z}\in\prod_{p\in\bP}\Qp$,
let $\mathbf{z}^{P}$ be the projection on $\prod_{p\in P}\Qp$ and
set according to the notation of Proposition \ref{pro:LLN-affQ} \[
\pi_{n}(\mathbf{z})=\pi_{n}^{P}(\mathbf{z}):=\pi_{n}(\mathbf{z}^{P}).\]
Fix an $\varepsilon>0$ and let $K=\sum_{p\not\in P}\pd^{-}+\varepsilon$.
Since $\nad{\mathbf{Z}_{\infty}^{*}}_{P}^{+}$ is finite,\[
\pr{x_{n}\in\mathcal{G}_{n\cdot K}^{\pi_{n}(\mathbf{Z}_{\infty}^{*})}}=\int_{B^{*}}\PP_{n}^{z}\left[\mathcal{G}_{n\cdot K}^{\pi_{n}(z)}\right]\nu^{*}(dz)\ver1\]
and $\PP_{n}^{z}\left[\mathcal{G}_{n\cdot K}^{\pi_{n}(z)}\right]$
converges to 1 for $\nu^{*}$-almost all $z$.

Let $h$ be the the $\PP^{z}$-almost sure limit of $-\ln\PP_{n}^{z}(x_{n})/n$
, which exists for $\nu^{*}$-almost all $z$ according to \cite{Ka00},
and consider the set\[
A_{n}=\left\{ g\in\affQ|-h-\varepsilon<\ln\PP_{n}^{z}(g)/n<-h+\varepsilon\right\} .\]
Then $\PP_{n}^{z}(A_{n}\cap\mathcal{G}_{n\cdot K}^{\pi_{n}(z)})$
converges to 1, while \[
\PP_{n}^{z}(A_{n}\cap\mathcal{G}_{n\cdot K}^{\pi_{n}(z)})\leq\e^{n\left(\varepsilon-h\right)}\card{\mathcal{G}_{n\cdot K}^{\pi_{n}(z)}}\leq\e^{n\left(\varepsilon-h\right)}\e^{(C+\varepsilon)\cdot n\cdot K}.\]
where $C$ is the parameter of the exponential growth of the gauges
$\mathcal{G}^{y}$. Thus, $(C+\varepsilon)\cdot K-h+\varepsilon\geq0$
and, since $\varepsilon$ was arbitrarily chosen, $h\leq C\cdot\sum_{p\not\in P}\pd^{-}.$
Now, we let $P$ grow to $P^{*}$ and we obtain \[
h\leq C\cdot\inf_{P\subseteq P^{*},\textrm{ finite }}\sum_{p\not\in P}\pd^{-}=C\cdot\sum_{p\not\in P^{*}}\pd^{-}=0.\]

\end{proof}
\bibliographystyle{alpha}
\bibliography{/home/brofferio/matematica/bound-rat/sara}

\end{document}